\documentclass[12pt]{amsart}

\usepackage{amssymb}
\usepackage[leqno]{amsmath}
\usepackage{pictex}

\newcommand\ZZ{\mathbb{Z}}
\newcommand\XX{\mathbb{X}}

\newcommand\GL{\mathbb{GL}}

\newcommand\diag{\mathrm{diag}}

\newcommand\rk{\operatorname{rk}}

\renewcommand\AA{\mathbb{A}}
\newcommand\DD{\mathbb{D}}
\newcommand\EE{\mathbb{E}}

\newcommand\transp{{\mathrm{tr}}}
\newcommand\REM[1]{}
\newcommand\Der{\operatorname{D^b}}
\newcommand\DER{{\mathcal D}}
\newcommand\md{\operatorname{mod}}
\renewcommand\mod{\operatorname{mod}}
\newcommand\coh{\operatorname{coh}}

\newcommand\Her{\mathcal{H}}
\newcommand\Clu{\mathcal{C}}
\newcommand\Tub{\mathcal{T}}
\newcommand\Cox{\Phi}

\newcommand\dual{\operatorname{D}}
\newcommand\Ext{\operatorname{Ext}}

\newcommand\Hom{\operatorname{Hom}}

\newcommand\Groth{\operatorname{K_0}}
\newcommand\oGroth{\operatorname{\overline{K}_0}}
\newcommand\spitz[1]{\langle #1\rangle}

\newcommand\gge[1]{\mathbf{#1}}
\newcommand\ogge[1]{\overline{\mathbf{#1}}}
\newcommand\cgge[1]{\widehat{\mathbf{#1}}}
\newcommand\lcm{\operatorname{lcm}}

\newcommand\TS{\mathcal S}
\newcommand\Sus{{\mathrm{T}}}
\newcommand\vollperiode[1]{v(#1)}
\newcommand\halbperiode[1]{h(#1)}

\newcommand\Dim[1]{|#1|}
\newcommand\Dimtwo[1]{|#1|_2}
\newcommand{\la}{\lambda}

\newcommand{\eeuler}[2]{\langle\langle #1,#2\rangle\rangle}
\newcommand{\euler}[2]{\langle#1,#2\rangle}

\newcommand{\eeulertwo}[2]{\eeuler{#1}{#2}_2}
\newcommand{\de}{\delta}
\newcommand{\De}{\Delta}
\newcommand{\bmu}{\overline{\mu}}
\newcommand{\bE}{\overline{E}}

\DeclareMathOperator\Coker{Coker}
\DeclareMathOperator\Imag{Im}
\DeclareMathOperator\Ker{Ker}

\newtheorem{Lemma}{Lemma}[section]
\newtheorem{Proposition}[Lemma]{Proposition}
\newtheorem{Theorem}[Lemma]{Theorem}

\theoremstyle{definition}

\newtheorem{nr}[Lemma]{\!\!}

\newcommand{\proofend}{\hfill$\Box$\par}

\newcommand{\HVCenter}[1]{\setbox 0=\hbox{#1}%
        \dimen0=\wd0%
        \dimen1=\ht0%
        \divide\dimen0 by 2%
        \divide\dimen1 by 2%
        \hskip -\dimen0%
        \lower \dimen1%
        \box0%
        \hskip -\dimen0}
\newcommand{\HBCenter}[1]{\setbox 0=\hbox{#1}%
        \dimen0=\wd0%
        \dimen1=\ht0%
        \divide\dimen0 by 2%
        \hskip -\dimen0%
        \box0%
        \hskip -\dimen0}
\newcommand{\HTCenter}[1]{\setbox 0=\hbox{#1}%
        \dimen0=\wd0%
        \dimen1=\ht0%
        \divide\dimen0 by 2%
        \hskip -\dimen0%
        \lower \dimen1%
        \box0%
        \hskip -\dimen0}
\newcommand{\RVCenter}[1]{\setbox 0=\hbox{#1}%
        \dimen0=\wd0%
        \dimen1=\ht0%
        \divide\dimen1 by 2%
        \hskip -\dimen0%
        \lower \dimen1%
        \box0%
        \hskip -\dimen0}
\newcommand{\RBCenter}[1]{\setbox 0=\hbox{#1}%
        \dimen0=\wd0%
        \hskip -\dimen0%
        \box0%
        \hskip -\dimen0}
\newcommand{\LVCenter}[1]{\setbox 0=\hbox{#1}%
        \dimen1=\ht0%
        \divide\dimen1 by 2%
        \lower \dimen1%
        \box0%
        \hskip -\dimen0}
\numberwithin{equation}{section}

\begin{document}
\parindent=0pt\parskip=4pt

  \begin{abstract}
    For the cluster category of a hereditary or a canonical
    algebra, equivalently for the cluster category of the hereditary category of
    coherent sheaves on a weighted projective line,
    we study the Grothendieck group with respect to an admissible
    triangulated structure.
  \end{abstract}

\title[Grothendieck group of a cluster category]{The Grothendieck
  group of a cluster category}
\author{M.~Barot}
\address{Instituto de Matem\'aticas\\
  Universidad Nacional Aut\'onoma de M\'exico\\
  Ciudad Universitaria, C.P. 04510\\
  MEXICO}
\email{barot\@@matem.unam.mx}
\author{D.~Kussin}
\address{Institut f\"ur Mathematik\\
Universit\"at Paderborn\\
33095 Paderborn\\
GERMANY}
\email{dirk@math.uni-paderborn.de}
\author{H.~Lenzing}
\address{Institut f\"ur Mathematik\\
Universit\"at Paderborn\\
33095 Paderborn\\
GERMANY}
\email{helmut@math.uni-paderborn.de}
\maketitle

\section{Introduction}

The cluster category $\Clu=\Clu (A)$ of a finite dimensional
hereditary algebra $A$ was introduced by Buan, Marsh, Reineke,
Reiten and Todorov~\cite{5clu}, in order to realize the cluster
algebras of Fomin and Zelevinsky~\cite{FZ1} via tilting theory.

The construction of the orbit category $\Clu(A)$, see \cite{Keller},
generalizes to the situation where $A$ is any
$k$-algebra of finite global dimension. In this paper, \emph{all}
algebras will be unitary, associative and of finite dimension over an
algebraically closed ground field $k$.

We call a triangulated structure $\TS$ on $\Clu$ \emph{admissible}
if the canonical projection functor $\pi:\Der(\md A)\rightarrow
\Clu$ is exact, that is, sends exact triangles to triangles from
$\TS$. We use the notation $\Clu_\TS$ if we consider $\Clu$ as a
triangulated category with triangulated structure $\TS$. We suspect
that an admissible triangulated structure for $\Clu$ may not be
unique.

By Keller \cite{Keller}, $\Clu$ admits an admissible triangulated
structure in case $\Der(\md A)$ is triangle-equivalent to
$\Der(\Her)$ for some hereditary abelian $k$-cate\-go\-ry $\Her$.
Assuming $\Her$ connected, by Happel's classification theorem this
happens if and only if $A$ is derived equivalent to a hereditary or
a canonical algebra, see \cite{Happel,HaRe}. In the first case, we
can choose $\Her=\md A$ where $A$ is hereditary and in the second
$\Her=\coh \XX$, the category of coherent sheaves over a weighted
projective line $\XX$, see \cite{GeLe}. In the present paper we
focus on the case $\Her=\coh\XX$, but also deal with the cases
$\Her=\md A$ where $A$ is the path algebra of a Dynkin or an
extended Dynkin quiver.

Given an admissible triangulated structure $\TS$ on $\Clu$ we study
the Gro\-then\-dieck group
$\Groth(\Clu_\TS)$ with respect to all triangles in $\TS$ and compare
it with the Grothendieck group $\oGroth(\Clu)$ with respect to all
\emph{induced} triangles, that is, the images of exact triangles of
$\Der(\md A)$ under the projection $\pi$.

Assuming $A$ of finite global dimension, we denote by $\Cox$ the
Coxeter transformation on $\Groth(\Der(\md A))$, that is, the map induced
by the Auslander-Reiten translation $\tau$ of $\Der(\md A)$.
In Section \ref{sec:oGroth} we show the following result.

\begin{Proposition}\label{prop:oGroth=coker}
  If $A$ is an algebra of finite global dimension and $\Clu=\Clu(A)$
  then we have $\oGroth(\Clu)=\Coker(1+\Cox)$.
\end{Proposition}

Let $A$ be a hereditary algebra of finite representation type or a
canonical algebra. In both cases $\oGroth(\Clu)$ and
$\Groth(\Clu_\TS)$ are shown to be free either over $\ZZ$ or over
$\ZZ_2=\ZZ/2\ZZ$ (independently of the admissible triangulated
structure $\TS$). We define the \emph{dual Grothendieck groups}
$\oGroth(\Clu)^\ast$ and $\Groth(\Clu)^\ast$ as the respective
$\ZZ$- or $\ZZ_2$-dual. In Section \ref{sec:Abschneiden} we show our
first main result.

\begin{Theorem}\label{thm:top}
We have $\Groth(\Clu_{\mathcal{S}})=\oGroth(\Clu)$ in each of the
following three cases:
\begin{itemize}
\item[(i)] $A$ is canonical with weight sequence $(p_1,\ldots,p_t)$
  having at least one even weight.
\item[(ii)] $A$ is tubular,
\item[(iii)] $A$ is hereditary of finite representation type.
\end{itemize}
\end{Theorem}

The remaining canonical cases are covered by the
next result.

\begin{Theorem}\label{thm:bottom}
Assume $\Clu=\Clu(A)$ is the cluster category of a canonical algebra
$A$ with weight sequence $(p_1,\ldots,p_t)$, where all weights $p_i$
are odd. For any admissible triangulated structure $\TS$ on $\Clu$
the Grothendieck group $\Groth(\Clu_\TS)$ is a non-zero quotient of
$\oGroth(\Clu)=\ZZ_2\oplus\ZZ_2$.
\end{Theorem}

Accordingly, if $A$ is canonical (of any weight type), we have
$\Groth(\Clu_\TS)\neq 0$ and Proposition \ref{prop:oGroth_cases}
yields an explicit basis of $\oGroth(\Clu)$. Since each tame
hereditary algebra is derived equivalent to a canonical one, Theorem
\ref{thm:top} and \ref{thm:bottom} cover also the tame hereditary
situation. To prove the two theorems our main device is to provide a
categorification of suitable members of the dual Grothendieck group
$\oGroth(\Clu)^\ast$, that is, to realize them by additive functions
on $\Clu_\TS$ in categorical terms of $\Clu$ (in a sense defined at
the beginning of Section~\ref{sec:Abschneiden}).

B.~Keller informed the authors that his
student Y.~Palu proved $\Groth(\Clu_\TS)=\oGroth(\Clu)$ for the
admissible structure $\TS$ constructed in \cite{Keller}.

In the last section, we consider the cluster category $\Clu(\Tub)$ of an
``isolated'' tube $\Tub$. We show that there always exists an admissible
triangulated structure on $\Clu(\Tub)$ and determine its Grothendieck
group explicitly.

\section{Notations and definitions}
\label{sec:Def}
\subsection*{Definition of cluster categories}
We assume that $A$ is an algebra (we recall that this means a
unitary, associative algebra of finite dimension over
$k=\overline{k}$) of finite global dimension. We denote by $\md A$
the category of finitely generated (or equivalently
finite-dimensional) right $A$-modules and by $\DER=\Der(\md A)$ the
bounded derived category of $\md A$. Since $A$ has finite global
dimension, $\DER$ is a triangulated category, see \cite{Happel1},
and we denote by $\Sus$ its suspension functor. Moreover, $\DER$ has
Auslander-Reiten triangles and the Auslander-Reiten translation
$\tau$ is an auto-equivalence of $\DER$.

Denoting $F=\tau^{-1}\circ \Sus$, the cluster category
$\Clu=\Clu(A)$ is defined as the \emph{orbit category}
$\Clu(A)=\DER/F^\ZZ$, whose objects are the objects of $\DER$ and
whose morphism spaces are given by
$$
\Hom_{\Clu(A)}(X,Y)=\bigoplus_{i\in\ZZ}\Hom_\DER(X,F^iY),
$$
which are finite dimensional spaces if $A$ is derived equivalent
  to a hereditary or to a canonical algebra.
We denote by $\pi:\DER\rightarrow\Clu(A)$ the canonical
projection functor and write occasionally $\pi X$ rather than $X$ for
objects in $\Clu$ for emphasis.

\subsection*{Admissible triangulated structures}
We call a triangulated structure $\TS$ on $\Clu$ \emph{admissible}
if the projection $\pi$ is exact and denote by $\Clu_\TS$ the
category $\Clu$ equipped with $\TS$. Keller \cite{Keller} proves the
existence of an admissible triangulated structure for $\Clu(A)$ if
$\Der(\md A)$ is triangle equivalent to $\Der(\Her)$ for some
hereditary abelian $k$-category $\Her$. Then $\Her$ has a
tilting complex, hence by \cite[Theorem 1.7]{HaRe} a tilting object.
We may assume that $\Her$ is connected. Passing to a derived
equivalent hereditary category we may then assume by Happel's
theorem \cite{Happel} that $\Her=\md H$, where $H$ is a hereditary
algebra, or $\Her=\coh \XX$, where $\XX$ is a weighted projective
line \cite{GeLe}. In
the first case $A$ is derived equivalent to a hereditary, in the
second case to a canonical algebra, see paragraph ``Canonical
algebras'' below. Since $\Clu$ -- up to equivalence -- only depends
on $\Der(\md A)$, we can assume that $A$ itself is hereditary or
canonical. Often, we also
shall write $\Clu(\Her)$ instead of $\Clu(A)$ if
  $\Der(\Her)\simeq\Der(\md A)$.

\subsection*{Grothendieck groups}
Any $\Clu$ as above is equipped with the auto\-equi\-va\-lence
$\tau:\Clu\rightarrow\Clu$, induced by the Auslander-Reiten
translation of $\Der(\md A)$. A triangle $X\rightarrow Y\rightarrow Z
\rightarrow \tau X$ in $\Clu$ is called \emph{induced} if it is
-- up to isomorphism -- the image under $\pi$ of an exact triangle in
$\Der(\md A)$. Note that $\tau$ takes the role of a suspension functor
for $\Clu$, although the induced triangles usually will not define a
triangulated structure on $\Clu$.
We denote by $\oGroth(\Clu)$ the Grothendieck group of $\Clu$ with
respect to all induced triangles.

If $\TS$ is an admissible triangulated structure on $\Clu$ we denote
by $\Groth(\Clu_\TS)$ the Grothendieck group of $\Clu$ with respect to
all triangles from $\TS$. Since each induced triangle lies in $\TS$ we
get a natural epimorphism
$$
\oGroth(\Clu)\rightarrow \Groth(\Clu_\TS).
$$

\subsection*{Hereditary categories}
If $\Her$ is hereditary then the derived category admits a simple
description: the indecomposable objects of $\Der(\Her)$ are of the form
$\Sus^i X$ for $X\in\Her$ indecomposable and some $i\in\ZZ$.
The morphism spaces are given by
\begin{equation}
  \label{eq:FR}
  \Hom_{\Der(\Her)}(\Sus^i X,\Sus^j Y)=\Ext_\Her^{j-i}(X,Y),\text{ for }X,Y\in\Her.
\end{equation}
In case $\Her=\coh \XX$, $\tau$ is an autoequivalence on $\Her$ and
therefore $\Her$ is a fundamental region for the functor
$F$, that is, for each indecomposable object $X\in\DER$ there exists a
unique $Y\in\Her$
such that $X=F^iY$ and therefore, we can identify the objects of
$\Clu$ with the objects of $\Her$ up to isomorphism.

We recall that the category $\coh \XX$ has \emph{Serre duality},
that is, there exists an autoequivalence $\tau$ for which
$\Ext^1_\Her(X,Y)\simeq \dual\Hom_\Her(Y,\tau X)$ holds functorially
in $X$ and $Y$. Similarly the categories $\DER=\Der(\coh\XX)$ and
$\DER=\Der(\md A)$, for $A$ hereditary, have also Serre duality in
the sense that $\Hom_\DER(X,\Sus Y)\simeq \dual\Hom_\DER(Y,\tau X)$
holds functorially in $X$ and $Y$.

\subsection*{Canonical algebras}
Canonical algebras were introduced by C.\ M.\ Ringel in \cite{Ri} as algebras
$A=k Q/I$, where the quiver $Q$ is obtained by joining a source $1$
with a sink $n$ by $t\geq 2$ arms consisting of $p_1,\ldots, p_t$ arrows
respectively, all pointing from $1$ to $n$:
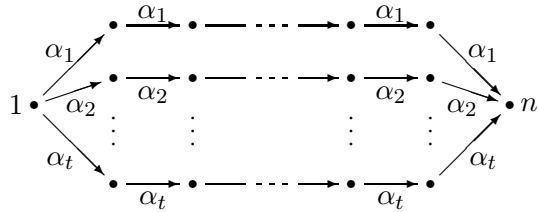
\begin{center}
  \begin{picture}(180,80)
    \put(0,10){
      \multiput(0,30)(180,0){2}{\circle*{3}}
      \multiput(0,60)(0,-20){2}{
        \multiput(30,0)(30,0){2}{\circle*{3}}
        \multiput(120,0)(30,0){2}{\circle*{3}}
        \multiput(35,0)(90,0){2}{\vector(1,0){20}}
        \put(65,0){\line(1,0){15}}
        \put(100,0){\vector(1,0){15}}
        \multiput(84,0)(5,0){3}{\line(1,0){2}}
        }
      \put(0,0){
        \multiput(30,0)(30,0){2}{\circle*{3}}
        \multiput(120,0)(30,0){2}{\circle*{3}}
        \multiput(35,0)(90,0){2}{\vector(1,0){20}}
        \put(65,0){\line(1,0){15}}
        \put(100,0){\vector(1,0){15}}
        \multiput(84,0)(5,0){3}{\line(1,0){2}}
        }
        \put(3.5,33.5){\vector(1,1){23}}
        \put(3.5,26.5){\vector(1,-1){23}}
        \put(153.5,3.5){\vector(1,1){23}}
        \put(153.5,56.5){\vector(1,-1){23}}
        \put(4.5,31.5){\vector(3,1){21}}
        \put(154.5,38.5){\vector(3,-1){21}}
        \put(10,48){\HBCenter{\small $\alpha_1$}}
        \multiput(45,63)(90,0){2}{\HBCenter{\small $\alpha_1$}}
        \put(170,48){\HBCenter{\small $\alpha_1$}}
        \put(18,32){\HTCenter{\small $\alpha_2$}}
        \multiput(45,37)(90,0){2}{\HTCenter{\small $\alpha_2$}}
        \put(162,32){\HTCenter{\small $\alpha_2$}}
        \put(10,12){\HTCenter{\small $\alpha_t$}}
        \multiput(45,-3)(90,0){2}{\HTCenter{\small $\alpha_t$}}
        \put(170,12){\HTCenter{\small $\alpha_t$}}
        \multiput(30,0)(30,0){2}{
          \multiput(0,15)(0,5){3}{\circle*{1}}
        }
        \multiput(120,0)(30,0){2}{
          \multiput(0,15)(0,5){3}{\circle*{1}}
        }
        \put(-4,30){\RVCenter{\small $1$}}
        \put(184,30){\LVCenter{\small $n$}}
      }
  \end{picture}
\end{center}
The ideal $I$ is generated by $t-2$ relations
$\alpha_i^{p_i}=\alpha_2^{p_2}-\mu_i \alpha_1^{p_1}$ for some pairwise
distinct $\mu_i\in k$ with $\mu_i\neq 0,1$.
The sequence $(p_1,\ldots,p_t)$ is called the \emph{weight sequence}
of $A$.
If $\sum_{i=1}^t\frac{1}{p_i}=t-2$ then $A$ is called
\emph{tubular}; this happens precisely for the weight sequences $(2,2,2,2)$,
$(3,3,3)$, $(2,4,4)$ and $(2,3,6)$.
We usually omit weights $p_i=1$ from the sequence, hence the weight
sequence $(3)$ means the sequence $(1,3)$.
We recall
that if $A$ is canonical of weight type $(p_1,\ldots,p_t)$ then
$\Der(\md A)\simeq \Der(\coh \XX)$ for a weighted projective
line $\XX$ of weight type $(p_1,\ldots,p_t)$.

\subsection*{Tubes}
  Let $\XX$ be a weighted projective line of weight type
  $(p_1,\ldots,p_t)$ and $\Her=\coh \XX$. We denote by $\Her_0$ the full
  subcategory of $\Her$ given by the objects of finite length and by
  $\Her_+$ the full subcategory of direct sums of indecomposable
  objects of infinite length.
  It is known, see \cite{GeLe}, that $\Her_0=\coprod_{x\in\XX}\Tub_x$
  is a coproduct of categories,
  where each $\Tub_x$ is a \emph{tube of rank} $q$, that is a
  connected, hereditary, uniserial category,
  which in abstract form can be realized as
  $\md_0^{\ZZ_q} k[[X]]$ (that is, as the category of $\ZZ_q$-graded
  $k[[X]]$-modules of finite length). Each tube of $\Her_0$ has rank one
  except finitely many (exceptional) tubes having rank $p_1,\ldots,p_t$,
  respectively.

  Furthermore $\Hom(\Her_0,\Her_+)=0$ and for each non-zero object
  $M\in\Her_+$ and each $x\in \XX$, we have $\Hom_\Her(M,\Tub_x)\neq 0$.

\subsection*{Formulas for $\Groth(\coh \XX)$}
The Grothendieck group $\Groth (\Her)$ of the abe\-lian category
$\Her=\coh \XX$
is described in detail in~\cite{K-Le,K-Ku}. It is equipped with the
Euler form defined by
$$
\spitz{[X],[Y]}=\dim_k \Hom_\Her (X,Y)-\dim_k \Ext^1_\Her (X,Y)
$$
on classes of objects $X$, $Y\in\Her$.
It follows from Serre duality that for all $\gge{x}$,
$\gge{y}\in\Groth (\Her)$ we have
$\spitz{\gge{y},\gge{x}}=-\spitz{\gge{x},\Cox\gge{y}}$, where $\Cox$
is the Coxeter transformation.

We denote by $L$ the structure sheaf and for each $i=1,\dots,t$
the unique simple sheaf $S_i$ belonging to the $i$-th
exceptional tube such that $\Hom_\Her (L,S_i)\neq 0$. Then
$\Hom_\Her (L,S_i)$ is one-dimensional, and $\Hom_\Her (L,\tau^j
S_i)=0$ for $j=1,\dots,p_i -1$. Furthermore, all simple sheaves from
homogeneous tubes have the same class in $\Groth(\Her)$;
we fix one, say $S_0$.
Now define the following elements of
$\Groth(\Her)$:
$$
\gge{a}=[L],\ \ \gge{s}_0=[S_0],\ \ \gge{s}_i=[S_i] \text{ for }
  i=1,\ldots,t.
$$
Define then the elements $\gge{s}_i(j)=\Cox^j\gge{s}_i$ for
$j\in\ZZ_{p_i}$. For later use we reproduce some facts from
\cite{K-Le}.
\begin{Proposition}\label{prop:Groth_Her}
Let $\Her=\coh\XX$ where $\XX$ is of weight type $(p_1,\ldots,p_t)$.
\begin{itemize}
\item[(a)]
  The abelian group $\Groth(\Her)$ is generated by the elements
  $\gge{a}$, $\gge{s}_0$, $\gge{s}_i(j)$,
  $i=1,\ldots,t$ and $j=0,\ldots,p_i-1$, subject to the defining
  relations
  \begin{equation}
    \label{eq:Groth0}
    \sum_{j=0}^{p_i-1}\gge{s}_i(j)=\gge{s}_0, \text{ for } i=1,\ldots,t.
  \end{equation}

\item[(b)]
  Define $p=\lcm(p_1,\ldots,p_t)$,
  $\delta=p\left(t-2-\sum_{i=1}^t
    \frac{1}{p_i}\right)$   and $\rk(\gge{x})=\spitz{\gge{x},\gge{s}_0}$.
Then for all $\gge{x}\in\Groth(\Her)$, we have

  \begin{equation*}
    \label{eq:Groth1}
  \Cox^p \gge{x}=\gge{x}+\delta\cdot \rk (\gge{x})\cdot \gge{s}_0
  \end{equation*}

\item[(c)] We have
  \begin{equation*}
    \label{eq:Groth2}
  \Cox\gge{a}=\gge{a}-\sum_{i=1}^t \gge{s}_i +(t-2)\cdot\gge{s}_0
  \end{equation*}
\item[(d)] Furthermore, we have
$\spitz{\gge{s}_i (m),\gge{s}_j (n)}=0$ for
$i\neq j$, and
$$
\spitz{\gge{s}_i (m),\gge{s}_i (n)}=
\begin{cases}
  1 & \text{ if } n\equiv m\mod p_i\\
  -1 & \text{ if } n\equiv m+1\mod p_i\\
  0 & \text{ else}
\end{cases}
$$
for all $i=1,\ldots,t$.\qed
\end{itemize}
\end{Proposition}

Beside the rank function $\rk(\gge{x})=\spitz{\gge{x},\gge{s}_0}$ we
also define the \emph{degree} function by
$$
\deg(\gge{x})=\sum_{j=0}^{p-1}\spitz{\Cox^j\gge{a},\gge{x}-\rk(\gge{x})\gge{a}}
$$
where $p=\lcm(p_1,\ldots,p_t)$. It is characterized by the properties
$\deg(L)=0$, $\deg(S_0)=p$ and $\deg(\tau^j S_i)=\frac{p}{p_i}$ for
$i=1,\ldots,t$ and $j\in\ZZ$.

\subsection*{Discriminant and slope}
\sloppy Let $\Her=\coh(\XX)$ be of weight type $(p_1,\ldots,p_t)$
and put $p=\lcm(p_1,\ldots,p_t)$. The \emph{discriminant}
$$
\de_\Her=p\bigl((t-2)-\sum_{i=1}^t 1/p_i\bigr)
$$
is an invariant of $\Her$ deciding on the complexity of the
classification problem for $\Her$, hence for $\Clu(\Her)$,
see~\cite{GeLe}. For $\de_\Her<0$ the category $\Her$ is derived
equivalent to the category $\md A$ for the path algebra $kQ$ of an
extended Dynkin quiver, and each such algebra $kQ$ has this
property. For $\de_\Her=0$ we are dealing with the tubular weights,
and for $\de_\Her>0$ the classification problem for $\Her$ is wild.
For this and the following statements we refer to \cite{GeLe}.

Each bundle $E$ has a \emph{line bundle filtration} $0=E_1\subset
E_1 \subset \cdots \subset E_r=E$ where each $E_i/E_{i-1}$ is a line
bundle. For each non-zero bundle $E$ its \emph{slope}
$\mu(E)=\deg(E)/\rk(E)$ is a rational number such that
$$
\mu(\tau E)=\mu(E)+\de_\Her
$$
holds. By means of line bundle filtrations for $E$ and $F$ it
follows that $\Hom_\Her(E,F)=0$ if $\mu(E)-\mu(F)$ is sufficiently
large. In particular, for $\de_\Her>0$ (resp.\ $\de_\Her<0$) we have
$\Hom_\Her(\tau^n E,F)=0$ (resp.\ $\Hom_\Her(E,\tau^nF)=0$ for
$n\gg0$. \fussy
\section{Grothendieck group with respect to induced triangles}
\label{sec:oGroth}

In this section, we describe the Grothendieck group $\oGroth(\Clu)$
with respect to the induced triangles.
Let $\DER=\Der(\md A)$. Then the Coxeter transformation
$\Cox:\Groth(\DER)\rightarrow\Groth(\DER)$ is given by
$\Cox([X])=[\tau X]$ for any object $X$ of $\DER$.

\subsection*{Proof of Proposition \ref{prop:oGroth=coker}}
The projection $\pi:\DER\rightarrow \Clu$ sends exact triangles to
induced triangles, hence yields an epimorphism
$$
\Groth(\DER)\rightarrow \oGroth(\Clu),[X]\mapsto [\pi X].
$$
We have $[F^{-1}X]=-[\tau X]$ in $\Groth(\DER)$, hence $[\pi X]=[\pi
F^{-1}X]=-[\pi \tau X]$ in $\oGroth(\Clu)$ showing that
$\pi(1+\Phi)=0$. In order to prove the exactness of
$$
\Groth(\DER)\xrightarrow{1+\Phi} \Groth(\DER)
\stackrel{\pi}\rightarrow \oGroth(\Clu)\rightarrow 0
$$
it therefore suffices to show that each morphism
$\lambda:\Groth(\DER)\rightarrow G$, for $G$ an abelian group, with
$\lambda(1+\Phi)=0$ induces a morphism
$\overline{\lambda}:\oGroth(\Clu)\rightarrow G$ with
$\lambda=\overline{\lambda} \pi$.

By the assumption $\lambda(1+\Phi)=0$ the corresponding function
$\lambda:\DER\rightarrow G$ is constant on $F$-orbits and additive
on exact triangles of $\DER$, hence induces a function
$\overline{\lambda}:\Clu\rightarrow G$ which is additive on induced
triangles. \qed

\subsection*{Explicit description of $\oGroth(\Clu)$}
Write $\ZZ_m$ for $\ZZ/m\ZZ$. We have the following general
description of $\oGroth(\Clu)$.

\begin{Proposition}\label{prop:explicit_oGroth}
Let $A$ be any algebra of finite global dimension and let $\Clu=\Clu(A)$.
Then $\oGroth(\Clu)$ has a unique expression as
$\ZZ^r\oplus\bigoplus_{i=1}^s (\ZZ_{m_i}\oplus
\ZZ_{m_i})$, for natural numbers $r,s$ and positive
$m_1,\ldots,m_s$ such that $m_i$ divides $m_{i+1}$ for all $i$.
Moreover, any such group occurs as $\oGroth(\Clu(H))$, where $H$ is
the hereditary path algebra given by the following quiver.
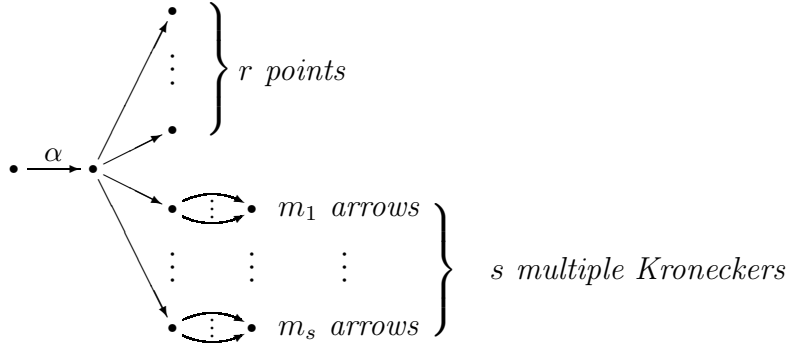
\begin{center}
\begin{picture}(280,120)
\put(0,60){
  \multiput(0,0)(30,0){2}{\circle*{3}}
  \put(15,3){\HBCenter{\small $\alpha$}}
  \put(5,0){\vector(1,0){20}}
  \multiput(60,15)(0,45){2}{\circle*{3}}
  \put(34,2){\vector(2,1){22}}
  \put(32,4){\vector(1,2){26}}
  \multiput(60,-15)(0,-45){2}{\circle*{3}}
  \multiput(90,-15)(0,-45){2}{\circle*{3}}
  \put(34,-2){\vector(2,-1){22}}
  \put(32,-4){\vector(1,-2){26}}
  \multiput(60,-15)(0,-45){2}{
    \qbezier(4,3)(15,8)(26,3)
    \put(27.3,2.2){\vector(2,-1){0.1}}
    \qbezier(4,-3)(15,-8)(26,-3)
    \put(27.3,-2.2){\vector(2,1){0.1}}
    \multiput(15,-3)(0,3){3}{\circle*{1}}
    }
  \multiput(60,32.5)(0,5){3}{\HVCenter{$\cdot$}}
  \multiput(60,-32.5)(0,-5){3}{\HVCenter{$\cdot$}}
  \multiput(90,-32.5)(0,-5){3}{\HVCenter{$\cdot$}}
  \multiput(125,-32.5)(0,-5){3}{\HVCenter{$\cdot$}}
  \put(70,35){$\left.\begin{picture}(1,27)\end{picture}\right\}$}
  \put(85,37.5){\LVCenter{$r$ points}}
  \put(100,-15){\LVCenter{$m_1$ arrows}}
  \put(100,-60){\LVCenter{$m_s$ arrows}}
  \put(155,-41){$\left.\begin{picture}(1,27)\end{picture}\right\}$}
  \put(180,-37.5){\LVCenter{$s$ multiple Kroneckers}}
}
\end{picture}
\end{center}
\end{Proposition}

Before we enter the proof we need some preparatory lemmas. If $Q$ is
a quiver, we denote by $B_Q$ the adjacency matrix of $Q$, that is,
$(B_Q)_{ij}$ denotes the number of arrows in $Q$ from $i$ to $j$.
Let $C$ be the Cartan matrix of $A$, that is, a matrix representing
the Euler form. Since $A$ has finite global dimension, $C$ has
determinant $\pm 1$ and $\Cox=-C^{-1}C^\transp$.

\begin{Lemma}\label{lem:adjacency}
  Let $A=k Q$ be a hereditary algebra and $\Clu=\Clu(A)$. Then we have
  $\oGroth(\Clu)=\Coker(B_Q-B_Q^\transp)$.
\end{Lemma}

\begin{proof}
  Since $A$ is finite-dimensional, $Q$ can not contain an oriented
  cycle. Hence the vertices of $Q$ can be ordered such that $C$ is
  upper triangular. Thus we see that $B_Q$ is nilpotent, hence
  $C=1+B_Q+B_Q^2+B_Q^3+\ldots$ is a finite sum and
  $C^{-1}=1-B_Q$. Therefore
  $1+\Cox=(C^{-\transp}-C^{-1})C^\transp=
  ((1-B_Q)^\transp-(1-B_Q))C^{\transp}=
  (B_Q-B_Q^\transp)C^\transp$,
  which shows that $\Coker (1+\Cox)= \Coker(B_Q-B_Q^\transp)$, thus
  the result follows by Proposition \ref{prop:oGroth=coker}.
\end{proof}

We call an arrow $\alpha:v\rightarrow w$ of a quiver $Q$ a \emph{source-arrow} if
$v$ is a source of $Q$ and $\alpha$ is the unique arrow of $Q$ starting in $v$.
Similarly an arrow $\alpha:w\rightarrow v$ is a \emph{sink-arrow} if
$v$ is a sink and $\alpha$ the unique arrow ending in $v$.
In both cases we denote by
$Q_{-\alpha}$ the quiver obtained from $Q$ by removing the
vertices $v$ and $w$ and all arrows starting or ending in $v$ or $w$.
The situation of a source-arrow is depicted as follows.
  
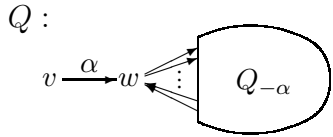
\begin{center}
    \begin{picture}(100,60)(0,-30)
      \put(0,0){\HVCenter{\small $v$}}
      \put(30,0){\HVCenter{\small $w$}}
      \put(5,0){\vector(1,0){20}}
      \put(15,3){\HBCenter{\small $\alpha$}}
      \put(36,2){\vector(2,1){20}}
      \put(36,1.33333){\vector(3,1){20}}
      \put(56,-12){\vector(-2,1){20}}
      \put(56,-8){\vector(-3,1){20}}
      \multiput(49,3)(0,-3){3}{\circle*{1}}
      \put(56,16){\line(0,-1){32}}
      \qbezier(56,16)(76,25)(96,16)
      \qbezier(56,-16)(76,-25)(96,-16)
      \qbezier(96,16)(106,11.5)(106,0)
      \qbezier(96,-16)(106,-11.5)(106,0)
      \put(81,0){\HVCenter{\small $Q_{-\alpha}$}}
      \put(0,20){\RBCenter{$Q:$}}
    \end{picture}
  \end{center}

The next result is quite useful for
calculating $\oGroth(\Clu)$ in practice.

\begin{Lemma}\label{lem:cutoff_arrow}
Let $Q$ be a quiver with an arrow $\alpha$, which is a source- or a
sink-arrow. Denote $H=k Q$ and $H'=k Q_{-\alpha}$. Then we have
$\oGroth(\Clu(H))\simeq \oGroth(\Clu(H'))$.
\end{Lemma}

\begin{proof}
  Assume that $\alpha$ is a source arrow (the case where $\alpha$ is a
  sink-arrow is similar).
  By renumbering the vertices, we can assume that $\alpha$ is the
  arrow $1\rightarrow 2$. Then we have
  $$
  B_Q-B_Q^\transp=\left[
    \begin{matrix}
      0 & 1 & 0 \\
      -1 & 0 & \rho^\transp \\
      0 & -\rho & B_{Q_{-\alpha}}-B_{Q_{-\alpha}}^\transp
    \end{matrix}
    \right].
    $$
    Adding multiples of the first row to the rows $3,\ldots,n$ and
    simultaneously adding (the same) multiples of the second column to
    the columns $3,\ldots,n$ we obtain a transformation matrix $T$ and
    a block diagonal matrix
    $$
    T(B_Q-B_Q^\transp)T^\transp=\diag(\left[
    \begin{matrix}
      0 & 1 \\
      -1 & 0
    \end{matrix}
    \right], B_{Q_{-\alpha}}-B_{Q_{-\alpha}}^\transp)
    $$
    and the result follows by Lemma \ref{lem:adjacency}.
\end{proof}

\subsection*{Proof of Proposition \ref{prop:explicit_oGroth}}
  If $C$ denotes the Cartan matrix of $A$ then
  $1+\Cox=(C^{-\transp}-C^{-1})C^\transp$. Now
  $S=(C^{-\transp}-C^{-1})$ is skew-symmetric. Clearly, we have
  $\Coker(1+\Cox)=\Coker S$.

  Using the skew-normal form of $S$, see \cite[Theorem IV.1]{Morris},
  we obtain  $S'=U^\transp S U$ for some $U\in\GL_n(\ZZ)$, where
  $S'=\diag(B_0,B_1,\ldots,B_s)$ is a
  block-diagonal matrix with the following blocks: $B_0$ is
  the zero matrix of size $r\times r$ and for $i=1,\ldots,s$,
  $$
  B_i=\left[
    \begin{matrix}
      0 & m_i \\ -m_i & 0
    \end{matrix}
  \right]
  $$
  where $m_i$ divides $m_{i+1}$ for all $i=1,\ldots,s-1$.  Therefore
  $\Imag S\simeq \Imag S'\simeq \bigoplus_{i=1}^s (m_i\ZZ)^2$ and we
  obtain $\oGroth(\Clu)\simeq \Coker S'\simeq \ZZ^r\oplus\bigoplus_{i=1}^s
  (\ZZ_{m_i}\oplus \ZZ_{m_i})$ as desired.

  Let $H$ be the hereditary algebra defined by the quiver in Proposition
  \ref{prop:explicit_oGroth} and denote $H'=k Q_{-\alpha}$. By Lemma
  \ref{lem:cutoff_arrow}, we have
  $\oGroth(\Clu(H))\simeq\oGroth(\Clu(H'))$. Now, the claim is obvious
  for $H'$ since
  $B_{Q_{-\alpha}}^{-1}-B_{Q_{\alpha}}^{-\transp}=\diag(B_0,B_1,\ldots,B_s)$ is the
  block-diagonal matrix as above.
\proofend

\subsection*{The hereditary case}

\begin{Proposition}
  If $A$ is a hereditary algebra whose quiver is a tree then
  $\oGroth(\Clu(A))$ is a free abelian group.
\end{Proposition}

\begin{proof}
  Any tree can be reduced to a disjoint union of $r$ vertices, for some
  $r$, by cutting off source- and sink-arrows.
  Hence, we get  $\oGroth(\Clu(A))\simeq\ZZ^r$ by Lemma \ref{lem:cutoff_arrow}.
\end{proof}

\begin{Proposition}\label{prop:hereditary2}
  Let $A$ be a connected hereditary representation-finite algebra,
  that is, the
  underlying graph of its quiver is a Dynkin diagram $\Delta$. Then, we have
  the following description.
  $$\oGroth(\Clu)=\begin{cases}
    0,&\ \ \text{ if $\Delta=\AA_n$, $\EE_n$ with $n$ even}\\
    \ZZ,&\ \ \text{ if $\Delta=\AA_n$, $\DD_n$, $\EE_7$ with $n$ odd}\\
    \ZZ^2,&\ \ \text{ if $\Delta=\DD_n$ with $n$ even}
    \end{cases}
  $$
\end{Proposition}

\begin{proof}
  This follows immediately using Lemma \ref{lem:cutoff_arrow}.
\end{proof}

\subsection*{The canonical case}
We now assume that $A$ is canonical of weight type
$(p_1,\ldots,p_t)$ and $\Her$ is the associated category of coherent
sheaves. We put $\Clu=\Clu(A)=\Clu(\Her)$ and start by describing $\oGroth(\Clu)$ by
generators and defining relations.

\begin{Proposition}
  The abelian group $\oGroth(\Clu)$ is generated by the elements
  $\ogge{a},\ogge{s}_0,\ogge{s}_1,\ldots\ogge{s}_t$ subject to the
  following defining relations.
  \begin{align}
  \label{eq:oGroth1} 2 \ogge{s}_0&=0,\\
  \label{eq:oGroth2} 2 \ogge{a}&=\sum_{i=1}^t (\ogge{s}_i-\ogge w),\\
  \label{eq:oGroth3} \ogge{s}_0& =\frac{1-(-1)^{p_i}}{2}\ogge{s}_i,
  \text{ for } i=1,\ldots,t.
  \end{align}
\end{Proposition}

\begin{proof}
  We recall from Proposition \ref{prop:Groth_Her}(a),
  that  $\Groth(\Her)$ is the abelian group generated by
  $\{\ogge{a},\ogge{s}_0,\ogge{s}_i(j)\mid i=1,\ldots,t\text{ and }
  j=0,\ldots,p_i-1\}$ subject to the defining relations \eqref{eq:Groth0}.
  Therefore $\oGroth(\Clu)=\Groth(\Her)/\Imag(1+\Cox)$ is the abelian
  group generated by the same generators with the relations
  \eqref{eq:Groth0} and the additional relations
  \begin{align}
    \label{eq:oGroth4} \ogge{a}+\Cox\ogge{a}&=0,\\
    \label{eq:oGroth5} \ogge{s}_0+\Cox\ogge{s}_0&=0\text { and }\\
    \label{eq:oGroth6} \ogge{s}_i(j)+\Cox\ogge{s}_i(j)&=0 \text{ for }
    i=1,\ldots,t\text{ and } j=1,\ldots,p_i,
  \end{align}
  which altogether form a system of defining relations.
  Using Proposition \ref{prop:Groth_Her}(c),
  we can rewrite \eqref{eq:oGroth4} as
  \eqref{eq:oGroth2}. Using $\Cox\gge{s}_0=\gge{s}_0$ we rewrite
  \eqref{eq:oGroth5} as \eqref{eq:oGroth1}.
  Using $\Cox\gge{s}_i(j)=\gge{s}_i(j+1)$ and \eqref{eq:Groth0} we
  obtain
  $$
  \ogge{s}_0=\sum_{j=0}^{p_1-1} (-1)^j \ogge{s}_i
  $$
  which can be rewritten in the form \eqref{eq:oGroth3}. Thus, since
  $\Cox\gge{s}_i(j)=\gge{s}_i(j+1)$, the group $\oGroth(\Clu)$ is
  generated by $\ogge{a}$, $\ogge{s}_0$, $\ogge{s}_i=\ogge{s}_i(0)$ for
  $i=1,\ldots,t$ subject to the defining relations \eqref{eq:oGroth1},
  \eqref{eq:oGroth2} and \eqref{eq:oGroth3}.
\end{proof}

\begin{Proposition}
  \label{prop:oGroth_cases}
  Let $\Her=\coh\XX$ with weight sequence $(p_1,\ldots,p_t)$ where
  $p_1,\ldots,p_r$ are even and $p_{r+1},\ldots,p_t$ are odd. Further let
  $\Clu=\Clu(\Her)$.
  \begin{itemize}
  \item[(i)] If $r\geq 1$ then
    $\oGroth(\Clu)$ is the free abelian group on
    $\ogge{a},\ogge{s}_2,\ldots,\ogge{s}_r$.
  \item[(ii)] If $r=0$ (that is, all weights $p_i$ are odd) then $\oGroth(\Clu)\simeq
    \ZZ \ogge{a}\oplus\ZZ\ogge{s}_0\simeq\ZZ_2\oplus\ZZ_2$.
  \end{itemize}
\end{Proposition}

\begin{proof}
  Let first $r\geq 1$. Then, by \eqref{eq:oGroth3}, we have
  $\ogge{s}_0=\frac{1-(-1)^{p_1}}{2}\ogge{s}_1=0$
  and for $i>r$, we obtain $\ogge{s}_i=0$, again by
  \eqref{eq:oGroth3}.
  Therefore $\ogge{s}_1=2\ogge{a}-\sum_{i=2}^r\ogge{s}_i$ because of
  \eqref{eq:oGroth2}. It follows that
  $\ogge{a},\ogge{s}_2,\ldots,\ogge{s}_r$ generate $\oGroth(\Clu)$
  without relations.

  Now let $r=0$. Then we obtain from \eqref{eq:oGroth3} that
  $\ogge{s}_i=\ogge{s}_0$ for all $i=1,\ldots,t$. Therefore we get that
  $\oGroth(\Clu)$ is generated by $\ogge{a}$ and $\ogge{s}_0$ with the
  remaining defining relations $2\ogge{s}_0=0$ and $2\ogge{a}=0$.
\end{proof}

\subsection*{The dual Grothendieck groups}
In the sequel the Grothendieck groups $\oGroth(\Clu)$ and
$\Groth(\Clu_\TS)$ are free over $\ZZ$ or $\ZZ_2$, respectively. We
define \emph{dual Grothendieck groups} $\oGroth(\Clu)^\ast$ and
$\Groth(\Clu_\TS)^\ast$ forming the respective $\ZZ$- or
$\ZZ_2$-duals.

We first deal with the $\ZZ$-free case. Since the Cartan
matrix has determinant $\pm 1$, the Euler form induces an
isomorphism $\Groth(\Her)\stackrel{\sim}\rightarrow\Groth
(\Her)^{\ast}$, $\gge{y}\mapsto\spitz{\gge{y},-}$. A linear form
$\lambda:\Groth (\Her)\rightarrow\ZZ$ induces a linear form
$\overline{\lambda}:\oGroth (\Clu)\rightarrow\ZZ$ if and only if
$\lambda\circ (1+\Cox)=0$.

\begin{Lemma}\label{lem:fix}
  A linear form $\lambda=\spitz{\gge{y},-}$ satisfies $\lambda\circ
  (1+\Cox)=0$ if and only if $\Phi\gge{y}=-\gge{y}$. In particular, in
  this case $\rk\gge{y}=0$ and $\deg\gge{y}=0$.
\end{Lemma}

\begin{proof}
  We have $\spitz{\gge{y},-}\circ (1+\Cox)=0$ if and only if
  $\spitz{\gge{y},\Cox^{-1}\gge{x}}+\spitz{\gge{y},\Cox\Cox^{-1}\gge{x}}=0$
  for all $\gge{x}\in\Groth (\Her)$, and since
  $\spitz{\gge{y},\Cox\gge{x}}=\spitz{\Cox^{-1}\gge{y},\gge{x}}$ this
  is equivalent to $\spitz{\gge{y}+\Cox\gge{y},-}=0$. Since
  the Cartan matrix has determinant $\pm 1$ the assertion follows.
\end{proof}

For any abelian group $G$ define $G_2=G\otimes_\ZZ \ZZ_2$.
Furthermore let $\rk_2,\deg_2:\Groth(\Her)_2\rightarrow \ZZ_2$ be
the functions induced by $\rk$ and $\deg$. Similarly define
$\spitz{-,-}_2:\Groth(\Her)_2\times\Groth(\Her)_2\rightarrow \ZZ_2$
to be induced by the Euler form.

\begin{Proposition}
  \label{prop:oGroth_stern}
  Let $\Her=\coh\XX$ with weight sequence $(p_1,\ldots,p_t)$ and set
  $\Clu=\Clu(\Her)$. The group $\spitz{\Cox}$ acts on $\Groth (\Her)$ by
  $\Cox .\gge{y}=-\Cox\gge{y}$.
  \begin{itemize}
  \item[(i)] If there is at least one even weight $p_i$ then
    there is an isomorphism $$\Groth
    (\Her)^{\spitz{\Cox}}\stackrel{\sim}\rightarrow\oGroth
    (\Clu)^{\ast},\ \gge{y}\mapsto\spitz{\gge{y},-}$$ which gives rise
    to an exact sequence $$0\rightarrow\oGroth
    (\Clu)^{\ast}\rightarrow\Groth
    (\Her)\stackrel{1+\Cox}\longrightarrow\Groth
    (\Her)\rightarrow\oGroth (\Clu) \rightarrow 0.$$
  \item[(ii)] If all weights are odd, then
  there is an isomorphism $$\Groth
    (\Her)_2^{\spitz{\Cox}}\stackrel{\sim}\rightarrow\oGroth
    (\Clu)^{\ast},\ \gge{y}\mapsto\spitz{\gge{y},-}_2$$ which gives rise
    to an exact sequence $$0\rightarrow\oGroth
    (\Clu)^{\ast}\rightarrow\Groth
    (\Her)_2\stackrel{1+\Cox}\longrightarrow\Groth
    (\Her)_2\rightarrow\oGroth (\Clu) \rightarrow 0.$$
  \end{itemize}
\end{Proposition}

\begin{proof}
  Part (i) follows from Lemma \ref{lem:fix} and the proof of (ii) is
  similar using reduction modulo $2$.
\end{proof}

If $\gge{x}\in\Groth(\Her)$ is a $\Cox$-periodic object with period
$q_{\gge{x}}$, we define
$$
\vollperiode{\gge{x}}=\sum_{j=0}^{q_{\gge{x}}-1}(-1)^j\Cox^j\gge{x}
$$
and if $q_{\gge{x}}$ is even, we define
$$
\halbperiode{\gge{x}}=\sum_{j=0}^{\frac{q_{\gge{x}}}{2}-1} \Cox^{2j}\gge{x}.
$$

\begin{Proposition}
  \label{prop:oGroth_stern_basis}
  Let $\Her=\coh\XX$ with weight sequence $(p_1,\ldots,p_t)$ where
  $p_1,\ldots,p_r$ are even and $p_{r+1},\ldots,p_t$ are odd. Further let
  $\Clu=\Clu(\Her)$.
  \begin{itemize}
  \item[(i)] If $r\geq 1$ then
  \begin{equation}
    \label{eq:basis_stern}
  \spitz{\vollperiode{\gge{s}_1},-},
  \spitz{\halbperiode{\gge{s}_2}-\halbperiode{\gge{s}_1},-},\ldots,
  \spitz{\halbperiode{\gge{s}_r}-\halbperiode{\gge{s}_1},-}
  \end{equation}
  is a $\ZZ$-basis of $\oGroth(\Clu)^\ast$.

  \item[(ii)] If $r=0$ (that is, all weights $p_i$ are odd) then
    $\rk_2$ and $\deg_2$ is a $\ZZ_2$-basis of $\oGroth(\Clu)^\ast$.
  \end{itemize}
\end{Proposition}

\begin{proof}
  (i)\ \
  Clearly $(1+\Cox)\vollperiode{\gge{s}_1}=0$ since $p_1$ is even.
  Furthermore, $(1+\Cox)\halbperiode{\gge{s}_i}=\gge{s}_0$ for
  $i=1,\ldots,r$ and therefore, by the Proposition
  \ref{prop:oGroth_stern}, we get that \eqref{eq:basis_stern} are
  indeed elements of $\oGroth(\Clu)^\ast$. From the formulas
  \begin{align*}
    \spitz{\vollperiode{\gge{s}_1},\gge{a}}&=1,
    &\spitz{\vollperiode{\gge{s}_1},\gge{s}_h}=&0\\
    \spitz{\halbperiode{\gge{s}_j}-\halbperiode{\gge{s}_1},\gge{a}}&=0,
    &\spitz{\halbperiode{\gge{s}_j}-\halbperiode{\gge{s}_1},\gge{s}_h}=&\delta_{jh}
  \end{align*}
  it follows that \eqref{eq:basis_stern} forms a $\ZZ$-basis of
  $\oGroth(\Clu)^\ast$.

  (ii)\ \
  We know that $\ogge{a},\ogge{s}_0$ is a $\ZZ_2$-basis of
  $\oGroth(\Clu)$ by Proposition \ref{prop:oGroth_cases}. We have
  $\Phi \gge{s}_0=\gge{s}_0$ and therefore $\rk_2=\spitz{-,\gge{s}_0}_2$
  defines a linear form on $\oGroth(\Her)$.

  Since $\ogge{s}_i=\ogge{s}_0$ for
  $i=1,\ldots,t$, we get from Proposition \ref{prop:Groth_Her}(d)
  that $\Phi \ogge{a}=\ogge{a} \mod 2$.
  Hence we get
  $$
  \deg_2(\gge{x})=\sum_{j=0}^p\spitz{\Phi^j\gge{a},\gge{x}-\rk(\gge{x})\gge{a}}_2
  = \spitz{\gge{a},\gge{x}-\rk(\gge{x})\gge{a}}_2=
  \spitz{\gge{a},\gge{x}}_2+\rk_2(\gge{x}).
  $$
  Thus, also $\deg$ induces a linear map
  $\deg_2:\oGroth(\Clu)\rightarrow\ZZ_2$.
  Since $\rk_2(\gge{a})= 1$,
  $\rk_2(\gge{s}_0)= 0$ and $\deg_2(\gge{a})= 0$,
  $\deg_2(\gge{s}_0)= 1$, it follows that $\rk_2$, $\deg_2$ form a
  $\ZZ_2$-basis of $\oGroth(\Clu)^\ast$.
\end{proof}

\section{Additive functions on $\Clu_\TS$}
\label{sec:Abschneiden}

\subsection*{Cutting technique}

For a finite dimensional $k$-vector space $V$ let $\Dim{V}$ (resp.\
$\Dimtwo{V}$) denote its $k$-dimension (resp.\ its $k$-dimension
modulo two). We put $\mu_E(X) = \Dim{\Hom_\Clu(E,X)}$ and write
$\bmu_E(X)$ for $\mu_E(X)$ modulo two.

In the sequel we identify members $\la$ from the dual Grothendieck
group $\oGroth(\Clu)^\ast$ with mappings $\la$ defined on
$\Clu=\Clu(\Her)$ with values in $\ZZ$, respectively in $\ZZ_2$,
that are additive on induced triangles. We call $\la$
\emph{realizable} if, depending on the case considered, it has the
form $\mu_E-\mu_F$ (resp.\ $\bmu_E$) with $E$ and $F$ from $\Clu$.
The realizable functions form a subgroup of $\oGroth(\Clu)^\ast$.
Note that usually neither $\mu_E$ nor $\bmu_E$ (respectively
$\mu_E-\mu_F$) are realizable. Our next proposition shows how to
construct realizable functions which additionally belong to
$\Groth(\Clu_\TS)^\ast$ for an admissible triangulated structure
$\TS$ on $\Clu$.

For any object $U\in\DER=\Der(\Her)$ and any positive integer
$q$ define the function $\lambda_X^{(q)}$ on the objects $Y$ of
$\Clu$ by
\begin{equation}
  \label{eq:lambda}
  \lambda_U^{(q)}:\Clu\rightarrow \ZZ,\
\lambda_U^{(q)}(Y)=\sum_{i=0}^{q-1}(-1)^i\Dim{\Hom_\Clu(\pi U,\Sus^i
Y)}
\end{equation}
and set $\overline{\lambda}_U^{(q)}:\Clu\rightarrow \ZZ_2,\ Y\mapsto
\lambda_U^{(q)}(Y) \mod 2$

\begin{Proposition}\label{prop:cutting}
  Suppose that $U$ is an object in $\Der(\Her)$ such that for some
  positive integer $q$ we have $\tau^q X\simeq F^m X$ for some
  $m\in\ZZ$.
\begin{itemize}
\item[(i)] If $q$ is even then
$\lambda_X^{(q)}$ is additive on each
triangle of an admissible triangulated structure on $\Clu$.
\item[(ii)] If $q$ is odd, then
$\overline{\lambda}_X^{(q)}$ is additive on each
triangle of an admissible triangulated structure on $\Clu$.
\end{itemize}
\end{Proposition}

\begin{proof}
  Identify $U$ with its image in $\Clu$. Let
  $X\stackrel{\alpha}\longrightarrow Y\stackrel{\beta}\longrightarrow
  Z\stackrel{\gamma}\longrightarrow \Sus X$ be a triangle in $\Clu$ with
  respect to an admissible triangulated structure. Application of the
  functor $\Hom_{\Clu}( U,-)$ gives a long exact sequence
  \begin{eqnarray*}
    \label{eq:long_exact:even}
    0\rightarrow K & \rightarrow &
    \Hom_\Clu(U,X)\rightarrow\Hom_\Clu(U,Y)
    \rightarrow\Hom_\Clu(U,Z)\rightarrow\\
    & \rightarrow & \Hom_\Clu(U,\tau X)\rightarrow
    \Hom_\Clu(U,\tau Y)\rightarrow\ \cdots\\
    & \cdots & \rightarrow\Hom_\Clu(U,\tau^{q-1} Y)\rightarrow
    \Hom_\Clu(U,\tau^{q-1}
    Z)\rightarrow K'\rightarrow 0,
\end{eqnarray*}
where $K=\Ker (\Clu (U,\alpha))$ and
\begin{equation}
  \label{eq:K=K'}
  K'= \Ker (\Clu (U,\tau^q
  \alpha))\simeq\Ker (\Clu (\tau^{-q} U,\alpha))\simeq\Ker (\Clu
  (U,\alpha))=K.
\end{equation}
The alternating sum of the dimensions of the spaces in the sequence
equals zero. If $q$ is even we hence get
\begin{equation}
  \label{eq:cU_additiv}
   \lambda_U^{(q)}(X) - \lambda_U^{(q)}(Y) + \lambda_U^{(q)}(Z) =0.
\end{equation}
Therefore $\lambda_U^{(q)}$ is a linear form on $\Groth (\Clu)$. If
$q$ is odd then this holds modulo $2$.
\end{proof}

\subsection*{The even canonical case}

We first study the case where $\Her=\coh \XX$ with weight sequence
$(p_1,\ldots,p_t)$. In this case we have linear forms which are
defined by ``periodic'' elements which lie in tubes: If
$U\in\mathcal{H}$ is indecomposable lying in a tube of rank $q$ then
$\tau^q U\simeq U$.

Assume that $p_1,\dots,p_r$ are even and $p_{r+1},\dots,p_t$ are
odd. Let $\mathcal{T}_1,\ldots,\mathcal{T}_r$ be the exceptional
tubes in $\Her_0$ of rank $p_1,\ldots,p_r$, respectively, and recall
that $S_i$ is a simple object from $\mathcal{T}_i$. By
Proposition \ref{prop:cutting} the functions
$\lambda_i=\lambda_{S_i}^{(p_i)}$ are additive on the triangles of
any admissible triangulated structure $\TS$ on $\Clu$.

If $\gge{x}$ is an element in $\Groth (\Her)$, denote by $\cgge{x}$
its image in $\Groth (\Clu_\TS)$.

\begin{Proposition}\label{prop:canonical_not_odd}
Assume that the number $r$ of even weights $p_i$ is non-zero, then
the linear forms $\lambda_i$ ($i=1,\ldots,r$) are realizable,
linearly independent over $\ZZ$ and
$\Groth(\Clu_\TS)=\oGroth(\Clu)\simeq \ZZ^r$.
\end{Proposition}

\begin{proof}
  Linear independence of $\lambda_1,\dots,\lambda_r$ follows from
  $\lambda_i (S_j)=2\delta_{ij}$, where $\delta_{ij}$ denotes the
  Kronecker symbol. We conclude that $\cgge{s}_1,\dots,\cgge{s}_r$ are
  linearly independent and hence also
  $\cgge{a},\,\cgge{s}_2,\dots,\cgge{s}_r$, since
  $2\cgge{a}=\sum_{i=1}^m \cgge{s}_i$. Therefore
  $\Groth(\Clu_\TS)=\oGroth(\Clu)\simeq \ZZ^r$ follows from
  Proposition \ref{prop:oGroth_cases} (i).
\end{proof}

\subsection*{The odd canonical case}

We adopt the notations of the previous section.  Recall that $S_0$
denotes a simple object from a homogeneous tube in $\Her_0$.

\begin{Proposition} \label{prop:realizable}\label{prop:canonical_odd}
Let $\Clu=\Clu(\Her)$, where $\Her$ is the category of coherent
sheaves on a weighted projective line of weight type
$(p_1,\ldots,p_t)$, where all weights $p_i$ are odd. Then the
following holds:

\begin{itemize}

\item[(i)] Always $\rk_2$ is a non-zero realizable member of
$\Groth(\Clu_\TS)^\ast\subseteq \oGroth(\Clu)^\ast$.

\item[(ii)] For $\de_\Her\neq0$ the subgroup of realizable members
of $\oGroth(\Clu)^\ast$ agrees with the subgroup $\langle
\rk_2\rangle$ generated by the rank modulo two.

\item[(iii)] For $\de_\Her=0$, that is for weight type $(3,3,3)$, we
have equality
$$
\Groth(\Clu_\TS)^\ast=
\oGroth(\Clu)^\ast=\ZZ_2\rk_2\oplus\ZZ_2\deg_2,
$$
and each member of $\oGroth(\Clu)^\ast$ is realizable.
\end{itemize}
\end{Proposition}

\begin{proof} (i) and (iii): We invoke Proposition~\ref{prop:oGroth_stern_basis}
and use that $\rk_2$ can be realized as
$\overline\lambda_0=\overline\lambda_{S_0}^{(1)}$ where
$\overline{\lambda}_0 (L)=1\ \mod 2$ and $\overline{\lambda}_0
(S_0)=0$.

In the tubular case, only the weight type $(3,3,3)$ matters, thus
the structure sheaf $L$ lies in a tube of $\tau$-period three. Hence
$\deg_2$ is realized by $\overline\lambda_L^{(3)}$ where
$\overline{\lambda}_L^{(3)} (S_0)=1\ \mod 2$.

(ii): Assume the function $\bmu_E(X)=|\Hom_\Clu(E,X)|_2$ with $E$
from $\Her$ is additive on induced triangles. We are going to show
that $\bmu_E$ is a multiple of $\rk_2$. Since $X\rightarrow
X\rightarrow 0\rightarrow \tau X$ is an induced triangle, we get
$$
\bmu_E(X)=\bmu_E(\tau X)=\bmu_{\tau^{-1} E}(X)
$$
for each object $X$ of $\Her$. Since $p=\lcm(p_1,\ldots,p_t)$ is
odd, setting $\bE=\bigoplus_{j=0}^{p-1}\tau^j E$ we thus obtain $
\bmu_E=\bmu_{\bE}$.

Next, we use a decomposition $E=E_+\oplus E_0$ of $E$ into a bundle
$E_+$ and an object $E_0$ of finite length. Invoking that $\tau^p$
acts as the identity on finite length objects of $\Her$, we see that
$\bE_0=\bigoplus_{j=0}^{p-1}\tau^j E$ is fixed under $\tau$. The
expression
$\bmu_{\bE_0}(X)=\Dimtwo{\Hom_\Her(\bE_0,X)}+\Dimtwo{\Ext^1_\Her(\bE_0,\tau^{-1}X)}$
hence agrees with $\euler{\bE_0}{X}_2$, and $\bmu_{\bE_0}$ is
a multiple of $\rk_2$. By part (i) the function $\rk_2$ is
additive on induced triangles, we conclude that the same holds for
$\bmu_{E_+}$. From now on, we may hence assume that $E$ is a bundle.
Note that
$$
\bmu_E(X)=\bmu_{\bE}(X)=\eeulertwo{E}{X} + \De_E(X),
$$
where $ \eeuler{E}{X}= \sum_{j=0}^{p-1}\euler{\tau^j E}{X}$,
$\eeulertwo{E}{X}=\eeuler{E}{X} \mod\ 2$ and
\begin{align*}
 \De_E(X)&=\sum_{j=0}^{p-1}\left(\Dimtwo{\Ext^1_\Her(\tau^j
E,X)}+\Dimtwo{\Ext^1_\Her(\tau^{j+1}E,X)}\right)\\
&=\Dimtwo{\Ext^1_\Her(E,X)}+\Dimtwo{\Ext^1_\Her(\tau^pE,X)}.
\end{align*}
By the Riemann-Roch formula,
 \begin{equation} \label{eq:RR}
\eeuler{E}{X}=-\frac{p}{2}\,\de_\Her\,\rk(E)\rk(X)+
\left|\begin{array}{cc}\rk(E)&\rk (X)\\ \deg (E)&\deg
(X)\end{array}\right|,
 \end{equation}

see \cite{K-Le}, the function $\eeulertwo{E}{-}$ is a linear
combination of $\rk_2$ and $\deg_2$. Hence $\eeulertwo{E}{-}$ is a
member of $\oGroth(\Clu)^\ast$, implying that $\De_E$ also belongs
to $\oGroth(\Clu)^\ast$.  By construction, $\De_E$ vanishes on
$S_0$. By means of a line bundle filtration of $E$, Serre duality
implies that $\De_E(L')=0$ for any line bundle $L'$ of sufficiently
large degree, and we deduce from Proposition~\ref{prop:oGroth_cases}
that $\De_E=0$. If $E$ is of even rank, then the function
$\eeulertwo{E}{-}$, hence also the function $\bmu_E$, is a multiple
of $\rk_2$, proving the claim in this case.

It remains to deal with the case that the rank of $E$ is
odd, where we deduce a contradiction from the assumption that
$\bmu_E$ belongs to $\oGroth(\Clu)^\ast$. Invoking $\bmu_E=\bmu_{\bE}$,
we have shown that $\De_{\bE}=0$. Note that $\De_{\bE}=0$ asserts
that
\begin{equation} \label{eq:ext}
\Dimtwo{\Ext^1_\Her(\bE,X)} = \Dimtwo{\Ext^1_\Her(\tau^{np}\bE,X)}
\end{equation}
for each $n\in \ZZ$ and each object $X$ from $\Her$.

\emph{Case $\de_\Her>0$}: Clearly, the functions
$\euler{\bE}{-}=\eeuler{E}{-}$ and
$\euler{\tau^{np}\bE}{-}=\eeuler{\tau^{np}E}{-}$ are additive on
induced triangles. They agree modulo two on $S_0$ and by formula
(\ref{eq:ext}) also on each line bundle $L'$ of large negative
degree. It then follows from Proposition~\ref{prop:oGroth_cases}
that $\euler{\bE}{-}_2=\euler{\tau^{np}\bE}{-}_2$.

By means of a line bundle filtration for $\bE$ we obtain for each
integer $n\gg 0$ two line bundles $L_1$ and $L_2$ of consecutive
degrees $d$ and $d+1$ such that
\begin{equation} \label{eq:contra}
\Hom_\Her (\tau^{np}\bE,L_i)=0 \textrm{ and } \Ext^1_\Her(\bE,L_i)=0
\textrm{ for } i=1,2.
\end{equation}
By \eqref{eq:RR} we get $\euler{\bE}{L_i}_2=\alpha + \deg_2(L_i)$
for some $\alpha\in\ZZ_2$, only depending on $E$. We then choose one
of the $L_i$ such that \eqref{eq:contra}  and further
$\euler{\bE}{L_i}_2=1$ holds. Invoking \eqref{eq:ext} we obtain the
contradiction
\begin{eqnarray*}
1&=&\euler{\bE}{L_i}_2=\euler{\tau^{np}\bE}{L_i}_2\\
&=&\Dimtwo{\Hom_\Her(\tau^{np}\bE,L_i)}+\Dimtwo{\Ext_\Her^1(\tau^{np}\bE,L_i)}=0.
\end{eqnarray*}

\emph{Case $\de_\Her<0$}: The proof is similar, choosing $n\ll0$.
\end{proof}

\subsection*{The Dynkin case}
Now, let $A$ be a connected hereditary repre\-sen\-ta\-tion-finite
algebra whose quiver has as underlying graph the Dyn\-kin diagram
$\Delta$. Then $\Delta$ is a star with length of the arms
$p_1,\ldots,p_t$ (where $t\leq 3$) and the Auslander-Reiten quiver
of $\Der(\md A)$ is $\ZZ\Delta$ whose $\tau$-orbits correspond to
the vertices of $\Delta$. In this case there are no tubes.
Nevertheless we find ``periodic'' objects. Let $m$ be the Coxeter
number of $\Delta$, that is, the order of the Coxeter transformation
$\Cox$. We have
$$m=
\begin{cases}
  n+1 & \text{ if }\Delta=\AA_n,\\
  2(n-1) & \text{ if }\Delta=\DD_n,\\
  12,\,18,\,30 & \text{ if }\Delta=\EE_6,\,\EE_7,\,\EE_8,\text{ respectively}.
\end{cases}$$

Since $\Groth (\Clu)=0$ in the cases $\Delta=\AA_n$ ($n$ even),
$\EE_6$, $\EE_8$ by Proposition~\ref{prop:hereditary2}, we restrict
our attention to the remaining cases.  Note that then $m$ is always an
even number.

\begin{Proposition}\label{prop:dynkin}
  In the cases $\Delta=\AA_n$ with $n$ odd or $\Delta=\EE_7$, let $M$
  be an indecomposable object of $\Der(\md A)$ lying in a $\tau$-orbit as
  indicated in the following picture.
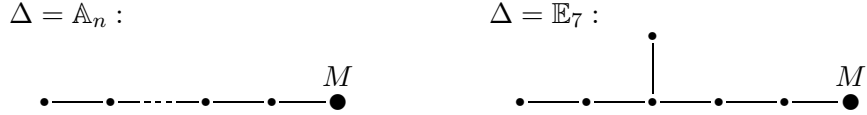
\begin{center}
\begin{picture}(300,45)
  \put(0,5){
    \multiput(0,0)(61,0){2}{
      \multiput(0,0)(25,0){2}{\circle*{3}}
      \put(3,0){\line(1,0){19}}
      }
    \multiput(28,0)(22,0){2}{\line(1,0){8}}
    \multiput(38,0)(4,0){3}{\line(1,0){2}}
    \put(89,0){\line(1,0){17}}
    \put(111,0){\circle*{6}}
    \put(111,6){\HBCenter{\small $M$}}
    \put(30,30){\RBCenter{\small $\Delta=\AA_n:$}}
    }
  \put(180,5){
    \multiput(0,0)(25,0){5}{\circle*{3}}
    \multiput(3,0)(25,0){4}{\line(1,0){19}}
    \put(50,3){\line(0,1){19}}
    \put(50,25){\circle*{3}}
    \put(125,0){\circle*{6}}
    \put(103,0){\line(1,0){17}}
    \put(125,6){\HBCenter{\small $M$}}
    \put(30,30){\RBCenter{\small $\Delta=\EE_7:$}}
    }
\end{picture}
\end{center}
Then $\lambda_M^{(m+2)}$ is a non-zero realizable function,
which is additive on all triangles.

In the case $\Delta=\DD_n$ ($n\geq 4$) one can choose indecomposable
objects $M_1$ and $M_2$ of $\Der(\md A)$ lying in the two $\tau$-orbits as
indicated in the following picture
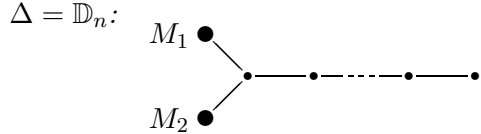
\begin{center}
  \begin{picture}(150,50)
    \put(90,20){
      \multiput(-16,-16)(0,32){2}{\circle*{6}}
      \multiput(0,0)(61,0){2}{
        \multiput(0,0)(25,0){2}{\circle*{3}}
        \put(3,0){\line(1,0){19}}
        }
      \multiput(28,0)(22,0){2}{\line(1,0){8}}
      \multiput(38,0)(4,0){3}{\line(1,0){2}}
      \put(-2,2){\line(-1,1){11}}
      \put(-2,-2){\line(-1,-1){11}}
      \put(-22,16){\RVCenter{\small{$M_1$}}}
      \put(-22,-16){\RVCenter{\small{$M_2$}}}
      }
    \put(0,40){\small $\Delta=\DD_n$:}
  \end{picture}
\end{center}
such that the functions $\lambda_i=\lambda_{M_i}^{(m+2)}$ for
$i=1,2$ are non-zero realizable and, for $n$ even, linearly
independent.
\end{Proposition}

\begin{proof}
  For any indecomposable $A$-module $U$ we have $\tau^m U\simeq \Sus^{-2}U$
  in $\DER=\Der(\md A)$ which implies $\tau^{m+2}U\simeq F^{-2}U$. By
  Proposition~\ref{prop:cutting} the function $\lambda_U^{(m+2)}$ is
  additive on triangles in $\Clu$ with respect to any admissible
  triangulated structure on $\Clu$.

  Let $\mathcal{H}=\mod A$. For indecomposable objects $M$ and $N$ in
  $\DER$ we have (identifying them with their images in $\Clu$)
  \begin{eqnarray*}
    \lambda_M^{(m+2)}(N) & = & \sum_{i=0}^{m+1}(-1)^i
    \Dim{\Hom_{\Clu}(\tau^{-i}M,N)}\\
    & = & \sum_{j\in\ZZ}\sum_{i=0}^{m+1}(-1)^i
    \Dim{\Hom_{\DER}(M,\tau^{i-j}\Sus^j N)}.
  \end{eqnarray*}
  \sloppy
  Setting $\mu_j (M,N)=\sum_{i=0}^{m+1}(-1)^i  \Dim{\Hom_{\DER}(M,\tau^{i-j}\Sus^j N)}$
  we have $\mu_j (M,N)=0$ for $j<0$ and
  $M$, $N\in\mathcal{H}\cup\tau^{-}\mathcal{H}$.

  \fussy
  In the case $\Delta=\AA_n$ ($n$ odd) we get
  $\lambda_M^{(m+2)}(M)=2$, where $M$ is as indicated above. Indeed,
  $\mu_0 (M,M)=1=\mu_1 (M,M)$ and $\mu_j (M,M)=0$ for $j\geq 2$.

  In the case $\Delta=\EE_7$ we have $\lambda_M^{(m+2)}(M)=6$. Indeed,
  $\mu_0 (M,M)=1$, $\mu_1 (M,M)=3$, $\mu_2 (M,M)=2$ and $\mu_j
  (M,M)=0$ for $j\geq 3$.

  In the case $\Delta=\DD_n$ let $M_1$ and $M_2$ be in the AR quiver
  lying in the following slice.

  
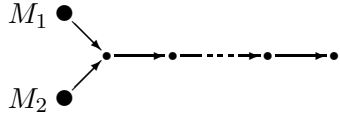
\begin{center}
  \begin{picture}(150,35)
    \put(90,20){
      \multiput(-16,-16)(0,32){2}{\circle*{6}}
      \multiput(0,0)(61,0){2}{
        \multiput(0,0)(25,0){2}{\circle*{3}}
        \put(3,0){\vector(1,0){19}}
        }
      \put(28,0){\line(1,0){8}}
      \put(50,0){\vector(1,0){8}}
      \multiput(38,0)(4,0){3}{\line(1,0){2}}
      \put(-13,13){\vector(1,-1){11}}
      \put(-13,-13){\vector(1,1){11}}
      \put(-22,16){\RVCenter{\small{$M_1$}}}
      \put(-22,-16){\RVCenter{\small{$M_2$}}}
      }
  \end{picture}
\end{center}

Let $\lambda_1 =\lambda_{M_1}^{(m+2)}$ and $\lambda_2
=\lambda_{M_2}^{(m+2)}$ where $M_1$ and $M_2$ are as indicated
above. Let $M$, $M' \in\{M_1,\,M_2\}$ with $M\neq M'$. It is easy to
see that $\Hom_{\DER}(M,\tau^i M)\neq 0$ if and only if $i$ is even
and $-(n-2)\leq i\leq 0$. Similarly, $\Hom_{\DER}(M,\tau^i M')\neq
0$ if and only if $i$ is odd and $1\leq i\leq n-1$. Moreover,
$$\tau^{-(n-1)}M\simeq
\begin{cases}
  \Sus M & n\ \text{even},\\
  \Sus M'& n\ \text{odd}.
\end{cases}$$ Using this we get $\mu_0 (M,M)=1$, $\mu_0 (M,M')=0$, and
$$\begin{array}{|c|c|c|c|c|}
  \hline
  & \mu_1 (M,M) & \mu_1 (M,M') & \mu_2 (M,M) & \mu_2 (M,M')\\
\hline
n\ \text{even} & \frac{n}{2} & -\frac{n-2}{2} & \frac{n-2}{2} &
-\frac{n-2}{2}\\
\hline
n\ \text{odd} & \frac{n-1}{2} & -\frac{n-1}{2} & \frac{n-3}{2} &
-\frac{n-1}{2}\\
\hline
\end{array}$$ and $\mu_j (M,M)=0=\mu_j (M,M')$ for $j\geq 3$.
Consequently, for even $n$ one has
$\lambda_1 (M_1)=n$, $\lambda_1 (M_2)=-(n-2)$, $\lambda_2 (M_1)=-(n-2)$ and
$\lambda_2 (M_2)=n$, and linear independence of $\lambda_1$ and
$\lambda_2$ follows. If $n$ is odd, then $\lambda_1
(M_1)=n-1=\lambda_2 (M_2)$ and $\lambda_1 (M_2)=-(n-1)=\lambda_2 (M_1)$.
\end{proof}

\subsection*{Proof of Theorem \ref{thm:top}}
For case (i) the assertion follows from Proposition
\ref{prop:canonical_not_odd}, for case (ii) it follows from
Proposition \ref{prop:canonical_odd} and for (iii) it follows from
the fact that by Proposition \ref{prop:oGroth=coker},
$\Groth(\Clu_\TS)$ is a quotient of $\oGroth(\Clu)$ and by
Proposition \ref{prop:dynkin} and \ref{prop:hereditary2} both are
free of the same rank.~\proofend

\subsection*{Proof of Theorem \ref{thm:bottom}}
This follows immediately from Propositions \ref{prop:canonical_not_odd} and
\ref{prop:canonical_odd}.
\proofend

\section{Cluster tubes}

\subsection*{Existence of admissible structures}
Let $\Tub$ be a tube of rank $q$. We consider the cluster category
$\Clu=\Clu(\Tub)$ as the orbit category $\Der(\Tub)/F^\ZZ$, where
again $F=\tau^{-1}\Sus$, where $\tau$ is the Auslander-Reiten
translation and $\Sus$ the suspension functor.  We call $\Clu(\Tub)$
the \emph{cluster tube} of rank $q$.

Since $\Tub$ has no tilting object, we can not invoke Keller's result
\cite{Keller} directly to conclude that $\Tub$ has an admissible
triangulated structure. We now show that $\Tub$ admits an admissible
structure anyway.

\begin{Proposition}
  The cluster tube $\Clu(\Tub)$ of rank $q$ admits an admissible
  triangulated structure.
\end{Proposition}

\begin{proof}
  Let $\XX$ be a weighted projective line of weight type
  $(q)=(1,q)$ and let $\Her=\coh \XX$. Recall the definitions of
  $\Her_0$ and $\Her_+$ from Section \ref{sec:Def}.
  We may view $\Tub$ as a full subcategory of $\Her_0$, which is even exact
  because of \eqref{eq:FR}.
  Therefore $\Clu(\Tub)$ is a full subcategory of $\Clu(\Her)$.

  By \cite{Keller}, there exists an admissible triangulated structure
  $\TS$ on $\Clu(\Her)$. We denote by $\TS'$ the subclass of $\TS$ given
  by all triangles $X\rightarrow Y\rightarrow Z\rightarrow \Sus X$ such
  that $X,Y,Z\in \Clu(\Tub)$. It is clear that once we show that
  $\TS'$ is a triangulated structure on $\Clu(\Tub)$ then it is
  admissible. Since $\Tub$, and then also $\Clu(\Tub)$,
  is closed under direct sums and summands in $\Her$, we only have to
  verify that $X,Y\in\Clu(\Tub)$ implies $Z\in\Clu(\Tub)$ for any
  triangle $X\rightarrow Y\rightarrow Z\rightarrow \Sus X$ in $\TS$.

  By the preceding remark, we can assume that $X,Y\in\Tub$ and
  $Z\in\Her$. Write $Z=Z_+\oplus Z_0$ where $Z_0\in\Her_0$ and
  $Z_+\in\Her_+$. Let $W\in\Her$ be a simple object in some
  homogeneous tube $\Tub'\neq \Tub$. Applying the functor
  $\Hom_{\Clu}(-,W)$ to the triangle $X\rightarrow Y\rightarrow
  Z\rightarrow \Sus X$, we get an exact sequence
  $$
  \Hom_\Clu(\Sus X,W)\rightarrow \Hom_\Clu(Z,W)\rightarrow
  \Hom_\Clu(Y,W)
  $$
  whose end terms are zero, because $\Tub$ and $\Tub'$ are orthogonal
  in $\Her$ and $\Clu(\Her)$. Therefore $\Hom_\Clu(Z,W)=0$, in
  particular $\Hom_\Her(Z,W)=0$. Hence $Z_+=0$ and $Z_0\not\in\Tub'$.
  Since we can vary
  $\Tub'\subset \Her_0$ we also see that $Z=Z_0\in\Tub$.
\end{proof}

\subsection*{The Grothendieck group of a cluster tube}

Let $\Tub$ be a tube and $\Clu=\Clu (\Tub)$ its cluster category. As
in~\ref{prop:oGroth=coker} one shows
$\oGroth(\Clu)=\Coker(1+\Cox)$.
We call an admissible triangulated structure on $\Tub$ an
\emph{induced} triangulated structure if it is obtained from an
embedding of $\Tub$ in $\coh \XX$ as explained in the previous paragraph.

\begin{Proposition}
  \label{prop:cluster_tube}
    Let $\Tub$ be a tube of rank $q$.
    \begin{itemize}
    \item[(i)] If $q$ is even then for any admissible triangulated
      structure $\TS$ on $\Clu=\Clu(\Tub)$ we have $\Groth(\Clu_\TS)=
      \oGroth(\Clu)\simeq\ZZ$.
    \item[(ii)]
      If $q$ is odd then
      for any induced triangulated structure
      $\TS$ on $\Clu=\Clu(\Tub)$ we have $\Groth(\Clu_\TS)=
      \oGroth(\Clu)\simeq\ZZ_2$.
    \end{itemize}
\end{Proposition}

\begin{proof}
  (i)\ \
  If $S$ is a simple object in $\Tub$ then $\Groth(\Tub)$ is the free
  group generated by the elements $\gge{s}(j)=[\tau^j S]$, for
  $j\in\ZZ_q$. Therefore $\ogge{s}(j)=-\ogge{s}(j+1)$ in
  $\oGroth(\Clu)$ and $\oGroth(\Clu)$ is generated by
  $\ogge{s}=\ogge{s}(0)$ without relation. This shows
  $\oGroth(\Clu)=\ZZ\ogge{s}\simeq \ZZ$.

  Finally, we can define $\lambda_S^{(q)}:\Clu\rightarrow \ZZ$ as in
  \eqref{eq:lambda} which defines a linear form
  $\lambda:\Groth(\Clu_\TS)\rightarrow \ZZ$ with $\lambda(S)=2$. Thus
  $\Groth(\Clu_\TS)$ has at least rank one and (i) follows.

  (ii)\ \
  Let $S\in\Tub$ be a simple object. Then $\oGroth (\Clu)$ is
  generated by $\ogge{s}$, where $\gge{s}=[S]$, and we have
  $2\ogge{s}=0$. We show that $\ogge{s}$ induces a non-trivial element
  in $\Groth (\Clu_\TS)$. For any object $X$ in $\Clu$ define
  $$\lambda (X)=\sum_{j=0}^{q-1}\Dimtwo{\Hom_{\Clu}(L,\tau^j X)}
  $$ For an object $X$ in $\Tub$ we have $\lambda (\pi X)=\deg_2
  (X)$. Indeed, since $\tau^q X\simeq X$
  \begin{eqnarray*}
    \lambda (\pi X) & = & \sum_{j=0}^{q-1}
    \Dimtwo{\Hom_{\Tub}(L,\tau^j X)}\pm\sum_{j=0}^{q-1}
    \Dimtwo{\Ext^1_{\Tub}(L,\tau^{j-1} X)}\\
     & = & \sum_{j=0}^{q-1}\spitz{L,\tau^j X}_2 =\deg_2 (X).
  \end{eqnarray*}
  In particular, $\lambda (\pi S)=1 \neq 0$. Now, $\lambda$ is
  additive on triangles in $\Clu$, which is shown with a version of
  the cutting technique similar to the proof of
  Proposition~\ref{prop:cutting}.
  In order to show that $K\simeq K'$ like in~\eqref{eq:K=K'} we use
  that $\tau^q$ is the identity functor on $\Tub$.
\end{proof}


\end{document}